\documentclass{amsart}

\usepackage{hyperref}
\usepackage{enumerate}
\usepackage{mathrsfs}
\usepackage{rotating}
\usepackage[matrix,arrow,graph,curve,frame]{xy}

\theoremstyle{plain}      
    \newtheorem{theorem}[section]{Theorem}
    \newtheorem{proposition}[section]{Proposition}
    \newtheorem{lemma}[section]{Lemma}
    \newtheorem{corollary}[section]{Corollary}
\theoremstyle{definition}
    \newtheorem{definition}[section]{Definition}
    \newtheorem{example}[section]{Example}
\theoremstyle{remark}

\newcommand{\A}{\ensuremath{\mathscr{A}}}
\newcommand{\B}{\ensuremath{\mathscr{B}}}
\newcommand{\C}{\ensuremath{\mathscr{C}}}

\newcommand{\Frob}{\ensuremath{\mathrm{Frob}}}

\newcommand{\ox}{\ensuremath{\otimes}}

\newcommand{\tint}{\ensuremath{\textstyle \int}}
\newcommand{\dint}{\ensuremath{\displaystyle \int}}

\newcommand{\ra}{\ensuremath{\xymatrix@=16pt@1{\ar[r]&}}}
\newcommand{\dra}{\ensuremath{\xymatrix@=20pt@1{\ar[r]&}}}

\begin{document}
%=========================================================================%
\title{Note on Frobenius monoidal functors}
\author{Brian Day}
\author{Craig Pastro}
\thanks{The first author gratefully acknowledges partial support of an
Australian Research Council grant while the second gratefully acknowledges
support of an international Macquarie University Research Scholarship and a
Scott Russell Johnson Memorial Scholarship.
The authors would like to thank Ross Street for several helpful comments.}

\address{Department of Mathematics \\
         Macquarie University \\
         New South Wales 2109 Australia}
\email{craig@ics.mq.edu.au}
\date{\today}

\begin{abstract}
It is well known that strong monoidal functors preserve duals. In this short
note we show that a slightly weaker version of functor, which we call
``Frobenius monoidal'', is sufficient.
\end{abstract}

\maketitle
%=========================================================================%

The idea of this note became apparent from Prop.~2.8 in the paper of
R.~Rosebrugh, N.~Sabadini, and R.F.C.~Walters~\cite{RSW}.

Throughout suppose that $\A$ and $\B$ are strict\footnote{We have decided
to work in the strict setting for simplicity of exposition, however, this
is not necessary.} monoidal categories.

\begin{definition}
A \emph{Frobenius monoidal functor} is a functor $F:\A \ra \B$ which is
monoidal $(F,r,r_0)$ and comonoidal $(F,i,i_0)$, and satisfies the
compatibility conditions
\begin{align*}
ir = (1 \ox r)(i \ox 1) &:F(A \ox B) \ox FC \dra FA \ox F(B \ox C) \\
ir = (r \ox 1)(1 \ox i) &:FA \ox F(B \ox C) \dra F(A \ox B) \ox FC,
\end{align*}
for all $A,B,C \in \A$.
\end{definition}

The compact case ($\ox$ = $\oplus$) of Cockett and Seely's linearly distributive
functors~\cite{CS} are precisely Frobenius monoidal functors, and Frobenius
monoidal functors with $ri = 1$ have been called \emph{split monoidal} by
Szlach\'anyi in~\cite{Szl}.

A \emph{dual situation} in $\A$ is a tuple $(A,B,e,n)$, where $A$ and $B$ are
objects of $\A$ and
\[
    e:A \ox B \dra I
    \qquad\qquad 
    n:I \dra B \ox A
\]
are morphisms in $\A$, called evaluation and coevaluation respectively,
satisfying the ``triangle identities'':
\[
    \xygraph{{A}="1"
        [r(2)] {A \ox B \ox A}="2"
        [d(1.2)] {A}="3"
        "1":"2" ^-{1 \ox n} 
        "2":"3" ^-{e \ox 1}
        "1":"3" _-1}
    \qquad\qquad
    \xygraph{{B}="1"
        [r(2)] {B \ox A \ox B}="2"
        [d(1.2)] {B.}="3"
        "1":"2" ^-{n \ox 1} 
        "2":"3" ^-{1 \ox e}
        "1":"3" _-1}
\]

\begin{theorem}\label{thm}
Frobenius monoidal functors preserve dual situations.
\end{theorem}

This theorem is actually a special case of the fact that linear functors
(between linear bicategories) preserve linear adjoints~\cite{CKS}.

\begin{proof}
Suppose that $(A,B,e,n)$ is dual situation in $\A$. We will show that
$(FA,FB,e,n)$, where $e$ and $n$ are defined as
\begin{align*}
    e &= \big(\xygraph{{FA \ox FB}
        :[r(2.1)] {F(A \ox B)} ^-r
        :[r(1.8)] {FI} ^-{Fe}
        :[r(1.1)] {I} ^-{i_0}}\big) \\
    n &= \big(\xygraph{{I}
        :[r(1.1)] {FI} ^-{r_0}
        :[r(1.8)] {F(B \ox A)} ^-{Fn}
        :[r(2.1)] {FB \ox FA} ^-{i}}\big),
\end{align*}
is a dual situation in $\B$.

The following diagram proves one of the triangle identities.
\[
    \xygraph{{FA}="1"
        [r(1.8)] {FA \ox FI}="2"
        [r(2.8)] {FA \ox F(B \ox A)}="3"
        [r(2.9)] {FA \ox FB \ox FA}="4"
        "2"[d(1.2)] {F(A \ox I)}="5"
        "3"[d(1.2)] {F(A \ox B \ox A)}="6"
        [d(1.2)] {F(I \ox A)}="8"
        "4"[d(1.2)] {F(A \ox B) \ox FA}="7"
        [d(1.2)] {FI \ox FA}="9"
        [d(1.2)] {FA}="10"
        "3"[d(0.6)r(1.5)] {\scriptstyle (\dagger)}
        "1":"2" ^-{1 \ox r_0}
        "2":"3" ^-{1 \ox Fn}
        "3":"4" ^-{1 \ox i}
        "5":"6" ^-{F(1 \ox n)}
        "6":"7" ^-i
        "8":"9" ^-i
        "1":"5" _-1
        "5":"8" _-1
        "8":"10" _-1
        "2":"5" ^-r
        "3":"6" ^-r
        "4":"7" ^-{r \ox 1}
        "6":"8" ^-{F(e \ox 1)}
        "7":"9" ^-{Fe \ox 1}
        "9":"10" ^-{i_0 \ox 1}}
\]
The square labelled by $(\dagger)$ requires the second Frobenius condition.
We remark that to prove the other triangle identity is similar and requires
the first Frobenius condition.
\end{proof}

\begin{proposition}
Any strong monoidal functor is a Frobenius monoidal functor.
\end{proposition}

\begin{proof}
Recall that a strong monoidal functor is a monoidal functor and a comonoidal
functor for which $r = i^{-1}$ and $r_0 = i_0^{-1}$. The commutativity of
the following diagram proves one of the Frobenius conditions.
\[
    \xygraph{{F(A \ox B) \ox FC}="1"
        [r(3.2)] {FA \ox FB \ox FC}="2"
        "1"[d(1.6)] {F(A \ox B \ox C)}="3"
        "2"[d(1.6)] {FA \ox F(B \ox C)}="4"
        "1":"2" ^-{i \ox 1}
        "3":"4" ^-i
        "1":@<-4pt>"3" _-r
        "3":@<-4pt>"1" _-i
        "2":@<4pt>"4" ^-{1 \ox r}
        "4":@<4pt>"2" ^-{1 \ox i}}
\]
The other is similar.
\end{proof}

\begin{proposition}
The composite of Frobenius monoidal functors is a Frobenius monoidal
functor.
\end{proposition}

\begin{proof}
Suppose that $F:\A \ra \B$ and $G:\B \ra \C$ are Frobenius monoidal
functors. It is well known and easy to see that the composite of monoidal
(resp. comonoidal) functors is monoidal (resp. comonoidal). We therefore
need only prove the Frobenius conditions, one of which follows from the
commutativity of
\[
    \xygraph{{GF(A \ox B) \ox GFC}="1"
        [r(3.5)] {G(F(A \ox B) \ox FC)}="2"
        [r(3.5)] {GF(A \ox B \ox C)}="3"
        "1"[d(1.2)] {G(FA \ox FB) \ox GFC}="4"
        [d(1.2)] {GFA \ox GFB \ox GFC}="7"
        "2"[d(1.2)] {G(FA \ox FB \ox FC)}="5"
        [d(1.2)] {GFA \ox G(FB \ox FC)}="8"
        "3"[d(1.2)] {G(FA \ox F(B \ox C))}="6"
        [d(1.2)] {GFA \ox GF(B \ox C),}="9"
        "2"[d(0.55)r(1.75)]{\scriptstyle (\ddagger)}
        "4"[d(0.55)r(1.75)]{\scriptstyle (\$)}
        "1":"2" ^-r
        "2":"3" ^-{Gr}
        "4":"5" ^-r
        "5":"6" ^-{G(1 \ox r)}
        "7":"8" ^-{1 \ox r}
        "8":"9" ^-{1 \ox Gr}
        "1":"4" _-{Gi \ox 1}
        "2":"5" _-{G(i \ox 1)}
        "4":"7" _-{i \ox 1}
        "5":"8" _-i
        "3":"6" ^-{Gi}
        "6":"9" ^-i}
\]
where the square labelled by $(\ddagger)$ uses the Frobenius property of $F$,
and the square labelled by $(\$)$ uses the Frobenius property of $G$.

The other Frobenius condition follows from a similar diagram.
\end{proof}

It is not too difficult to see that a Frobenius monoidal functor
$F:\mathbf{1} \ra \A$ is a Frobenius algebra in $\A$. Therefore, we have
the following corollary.

\begin{corollary}
Frobenius monoidal functors preserve Frobenius algebras. That is, if $R$ is a
Frobenius algebra in $\A$ and $F:\A \ra \B$ is a Frobenius functor, then
$FR$ is a Frobenius algebra in $\B$.
\end{corollary}

\begin{example}
Suppose that $\A$ is a braided monoidal category. If $R \in \A$ is a
Frobenius algebra in $\A$, then $F = R \ox -:\A \ra \A$ is a Frobenius
monoidal functor. The monoidal structure $(F,r,r_0)$ is given by
\begin{align*}
r_{A,B} &= \big(\xygraph{{R \ox A \ox R \ox B}
        :[r(3)] {R \ox R \ox A \ox B} ^-{1 \ox c \ox 1}
        :[r(2.8)] {R \ox A \ox B} ^-{\mu \ox 1 \ox 1}}\big)  \\
r_0 &= \big(\xygraph{{I}
        :[r(1)] {R} ^-\eta}\big)
\end{align*}
and the comonoidal structure $(F,i,i_0)$ by
\begin{align*}
i_{A,B} &= \big(\xygraph{{R \ox A \ox B}
        :[r(2.8)] {R \ox R \ox A \ox B} ^-{\delta \ox 1 \ox 1}
        :[r(3)] {R \ox A \ox R \ox B} ^-{1 \ox c \ox 1}}\big)  \\
i_0 &= \big(\xygraph{{R}
        :[r(1)] {I} ^-\epsilon}\big).
\end{align*}
The Frobenius conditions now follow easily from the properties of
Frobenius algebras.

This example shows that Frobenius monoidal functors generalize Frobenius
algebras much in the same way that monoidal comonads, or comonoidal monads,
generalize bialgebras.
\end{example}

The following proposition is a generalization of the fact that morphisms of
Frobenius algebras (morphisms which are both algebra and coalgebra morphisms)
are isomorphisms. It also generalizes the result that monoidal natural
transformations between strong monoidal functors with (left or right) compact
domain are invertible.

\begin{proposition}
Suppose that $F,G:\A \ra \B$ are Frobenius monoidal functors and that
$\alpha:F \ra G$ is a monoidal and comonoidal natural transformation. If
$A \in \A$ is part of a dual situation, i.e., $(A,B,e,n)$ or $(B,A,e,n)$
is a dual situation, then $\alpha_A:FA \ra GA$ is invertible.
\end{proposition}

\begin{proof}
We shall assume that $A$ is part of the dual situation $(A,B,e,n)$. The
other case is treated similarly. The component $\alpha_B:FB \ra GB$ has mate
\[
    \xygraph{{GA}
        :[r(2.2)] {GA \ox FB \ox FA} ^-{1 \ox n}
        :[r(3.4)] {GA \ox GB \ox FA} ^-{1 \ox \alpha_B \ox 1}
        :[r(2.2)] {FA} ^-{e \ox 1}}
\]
which we will show is the inverse to $\alpha_A$.

If $\alpha$ is both monoidal and comonoidal then the diagrams
\[
    \vcenter{\xygraph{{FA \ox FB}="1"
        [d(1.2)] {F(A \ox B)}="2"
        [d(1.2)] {FI}="3"
        [d(1)r(1.3)] {I}="x"
        "1"[r(2.6)] {GA \ox GB}="4"
        "2"[r(2.6)] {G(A \ox B)}="5"
        "3"[r(2.6)] {GI}="6"
        "1":"2" _-r
        "2":"3" _-{Fe}
        "4":"5" ^-r
        "5":"6" ^-{Ge}
        "1":"4" ^-{\alpha_A \ox \alpha_B}
        "2":"5" ^-{\alpha_{A \ox B}}
        "3":"6" ^-{\alpha_I}
        "3":"x" _-{i_0}
        "6":"x" ^-{i_0}}}
\qquad\qquad
    \vcenter{\xygraph{{FB \ox FA}="1"
        [u(1.2)] {F(B \ox A)}="2"
        [u(1.2)] {FI}="3"
        [u(1)r(1.3)] {I}="x"
        "1"[r(2.6)] {GB \ox GA}="4"
        "2"[r(2.6)] {G(B \ox A)}="5"
        "3"[r(2.6)] {GI}="6"
        "2":"1" _-i
        "3":"2" _-{Fn}
        "5":"4" ^-i
        "6":"5" ^-{Gn}
        "1":"4" ^-{\alpha_A \ox \alpha_B}
        "2":"5" ^-{\alpha_{A \ox B}}
        "3":"6" ^-{\alpha_I}
        "x":"3" _-{r_0}
        "x":"6" ^-{r_0}}}
\]
commute. The following diagrams prove that $\alpha_A$ is invertible.The
first diagram above says exactly that the triangle labelled by $(\pounds)$
below commutes. The second diagram above that the triangle labelled by
$(\yen)$ below commutes.
\[
    \xygraph{{FA}="1"
        [r(2.9)] {GA}="2"
        "1"[d(1.2)] {FA \ox FB \ox FA}="3"
        [d(1.2)] {FA}="5"
        "2"[d(1.2)] {GA \ox FB \ox FA}="4"
        [d(1.2)] {GA \ox GB \ox FA}="6"
        "3"[d(0.8)r(0.7)]{\scriptstyle (\pounds)}
        "1":"2" ^-\alpha
        "3":"4" ^-{\alpha \ox 1 \ox 1}
        "6":"5" ^-{e \ox 1}
        "1":"3" _-{1 \ox n}
        "3":"5" _-{e \ox 1}
        "2":"4" ^-{1 \ox n}
        "4":"6" ^-{1 \ox \alpha \ox 1}
        "3":"6" ^-{\alpha \ox \alpha \ox 1}}
\ \ 
    \xygraph{{GA}="1"
        [r(2.9)] {GA \ox FB \ox FA}="2"
        "1"[d(1.2)] {GA \ox GB \ox GA}="3"
        [d(1.2)] {GA}="5"
        "2"[d(1.2)] {GA \ox GB \ox FA}="4"
        [d(1.2)] {FA}="6"
        "1"[d(0.4)r(0.7)]{\scriptstyle (\yen)}
        "1":"2" ^-{1 \ox n}
        "4":"3" ^-{1 \ox 1 \ox \alpha}
        "6":"5" _-\alpha
        "2":"3" ^-{1 \ox \alpha \ox \alpha}
        "1":"3" _-{1 \ox n}
        "3":"5" _-{e \ox 1}
        "2":"4" ^-{1 \ox \alpha \ox 1}
        "4":"6" ^-{e \ox 1}}
\]
\end{proof}

Denote by $\Frob(\A,\B)$ the category of Frobenius monoidal functors from
$\A$ to $\B$ and all natural transformations between them.

\begin{proposition}[cf.~\cite{RSW} Prop.~2.10]
If $\B$ is a braided monidal category, then $\Frob(\A,\B)$ is a braided
monoidal category with the pointwise tensor product of functors.
\end{proposition}

\begin{proof}
Consider the pointwise tensor product of Frobenius monoidal functors
$F,G:\A \ra \B$. That is,
\[
    (F \ox G)A = FA \ox GA.
\]
It is obviously an associative and unital tensor product with unit $I(A) = I$
for all $A \in \A$.

We may define morphisms as follows:
\begin{align*}
    r &= (r \ox r)(1 \ox c^{-1} \ox 1)
        :(F \ox G)A \ox (F \ox G)B \dra (F \ox G)(A \ox B) \\
    r_0 &= r_0 \ox r_0:I \dra (F \ox G)I \\
    i &= (1 \ox c \ox 1)(i \ox i)
        :(F \ox G)(A \ox B) \dra (F \ox G)A \ox (F \ox G)B \\
    i_0 &= i_0 \ox i_0:(F \ox G)I \dra I.
\end{align*}
That these morphisms provide a monoidal and a comonoidal structure on
$F \ox G$ is not too difficult to show, and is omitted here. The following
diagram proves the first Frobenius condition, where the $\ox$ symbol has
been removed as a space spacing mechanism.
\[
    \scalebox{0.85}{
    \xygraph{{F(AB) ~ G(AB) ~ FC ~ GC}="1"
        [r(4)] {FA ~ FB ~ GA ~ GB ~ FC ~ GC}="2"
        [r(4.7)] {FA ~ GA ~ FB ~ GB ~ FC ~ GC}="3"
        "1"[d(1.2)] {F(AB) ~ FC ~ G(AB) ~ GC}="4"
        [d(1.2)] {F(ABC) ~ G(ABC)}="7"
        "2"[d(1.2)] {FA ~ FB ~ FC ~ GA ~ GB ~ GC}="5"
        [d(1.2)] {FA ~ F(BC) ~ GA ~ G(BC)}="8"
        "3"[d(1.2)] {FA ~ GA ~ FB ~ FC ~ GB ~ GC}="6"
        [d(1.2)] {FA ~ GA ~ F(BC) ~ G(BC)}="9"
        "1":@{=}[u(0.7)] {(F\! \ox\! G)(AB)\! \ox\! (F\! \ox\! G)C}
        "9":@{=}[d(0.7)] {(F\! \ox\! G)A\! \ox\! (F\! \ox\! G)(BC)}
        "1":"2" ^-{i \> i \> 1 \> 1}
        "2":"3" ^-{1 \> c \> 1 \> 1 \> 1}
        "4":"5" ^-{i \> 1 \> i \> 1}
        "5":"6" ^-{1 \> c_{FBFC,GA} \> 1 \> 1}
        "7":"8" ^-{i \> i}
        "8":"9" ^-{1 \> c \> 1}
        "1":"4" _-{1 \> c^{-1} \> 1}
        "2":"5" _-{1 \> 1 \> c_{GAGB,FC}^{-1} \> 1}
        "3":"6" ^-{1 \> 1 \> 1 \> c^{-1} \> 1}
        "4":"7" _-{r \> r}
        "5":"8" _-{1 \> r \> 1 \> r}
        "6":"9" ^-{1 \> 1 \> r \> r}}}
\]
The bottom left square commutes by the Frobenius condition, and the others
by properties of the braiding. The second Frobenius condition follows from
a similar diagram. So, $F \ox G$ is a Frobenius monoidal functor.

The braiding $c_{F,G}:F \ox G \ra G \ox F$ is given on components by
\[
    (c_{F,G})_A = c_{FA,GA}:FA \ox GA \dra GA \ox FA.
\]
\end{proof}

\begin{corollary}
If $\B$ is a braided monoidal category and $\A$ is a self-dual compact
category, meaning that for any object $A \in \A$, $(A,A,e,n)$ is a dual
situation in $\A$, then $\Frob(\A,\B)$ is a self-dual braided compact category.
\end{corollary}

\begin{proof}
By Theorem~\ref{thm} Frobenius monoidal functors preserve duals, and
therefore, for any $A \in \A$, $(FA,FA,e,n)$ is a dual situation in $\B$.
\end{proof}

Recall that, if $\A$ is a small monoidal category, and if small colimits
exist and commute with the tensor product in $\B$, then the equations
\begin{align*}
    F*G &= \int^{A,B} \A(A \ox B,-) \cdot FA \ox FB \\
    J &= \A(I,-) \cdot I,
\end{align*}
where $\cdot$ denotes copower, describe the \emph{convolution monoidal
structure} on the functor category $[\A,\B]$ (cf.~\cite{D}). Then we have:

\begin{theorem}\label{thm-fun}
If $\A$ is a small monoidal category and $\B$ is a
monoidal category having all small colimits commuting with tensor, then any
Frobenius monoidal functor $F:\A \ra \B$ for which the canonical evaluation
morphism
\begin{equation}\label{3to1}\tag{$\flat$}
\int^{A,B,C} \A(A \ox B \ox C,-) \cdot F(A \ox B \ox C) \dra F
\end{equation}
is an isomorphism, becomes an algebra with a comultiplication which satisfies
the Frobenius identities in the convolution functor category $[\A,\B]$.

Note that, by the Yoneda lemma, the equation~\eqref{3to1} is satisfied by
all the functors $F:\A \ra \B$ if $\A$ is a closed monoidal category and
the canonical evaluation morphism
\begin{equation}\label{eq-asump}\tag{$\sharp$}
    \int^{B,C} \A(A,B \ox C \ox [B \ox C,-]) \dra \A(A,-)
\end{equation}
is an isomorphism for all $A \in \A$.
\end{theorem}

Before we prove Theorem~\ref{thm-fun} we will need the following lemma.

\begin{lemma}\label{lem-2to1}
Assuming equation~\eqref{3to1} in Theorem~\ref{thm-fun}, we may also derive
the two variable version, that is, that the canonical evaluation morphism
\[
    \int^{A,B} \A(A \ox B,-) \cdot F(A \ox B) \dra F
\]
is an isomorphism.
\end{lemma}

\begin{proof}
The canonical evaluation morphism
\[
    \xygraph{{\dint^{A,B,C} \A(A \ox B \ox C,-) \cdot F(A \ox B \ox C)}
        :[r(5.3)] {\dint^{A,B} \A(A \ox B,-) \cdot F(A \ox B)} ^-h}
\]
is a retraction of (either of the canonical morphisms in the opposite
direction), say, $k$.
We may compose the canonical morphism
\[
    \int^{A,B} \A(A \ox B,-) \cdot F(A \ox B) \dra F
\]
with the isomorphism
\[
    \xygraph{{F}
    :[r(3.6)] {\dint^{A,B,C} \A(A \ox B \ox C,-) \cdot F(A \ox B \ox C)} ^-\cong}
\]
of our assumption to get a morphism
\[
    \xygraph{{\dint^{A,B} \A(A \ox B,-) \cdot F(A \ox B)}
    :[r(5.3)] {\dint^{A,B,C} \A(A \ox B \ox C,-) \cdot F(A \ox B \ox C)} ^-l},
\]
which makes the diagram
\[
    \xygraph{{\A(A \ox B \ox C,-) \cdot F(A \ox B \ox C)}="x"
      "x"[r(5)] {\dint^{A,B} \A(A \ox B,-) \cdot F(A \ox B)}="2"
      "2"[u(1.6)] {\dint^{A,B,C} \A(A \ox B \ox C,-) \cdot F(A \ox B \ox C)}="1"
      "2"[d(1.6)] {\dint^{A,B,C} \A(A \ox B \ox C,-) \cdot F(A \ox B \ox C)}="3"
        "x":"1" ^-{\mathrm{copr}}
        "x":"2" ^-{\mathrm{copr}}
        "x":"3" _-{\mathrm{copr}}
        "1":"2" ^-h
        "2":"3" ^-l}
\]
commute. Therefore $lh=1$. We have
\[
    hl = hlhk = hk = 1
\]
so $l$ is an isomorphism, hence the canonical evaluation morphism
\[
    \int^{A,B} \A(A \ox B,-) \cdot F(A \ox B) \dra F
\]
is an isomorphism.
\end{proof}

A consequence of Lemma~\ref{lem-2to1} is that we may write
\begin{align*}
    F * F &=     \int^{X,C} \A(X \ox C,-) \cdot FX \ox FC \\
          &\cong \int^{X,C} \A(X \ox C,-) \cdot \Big(\int^{A,B} \A(A \ox B,X)
                    \cdot F(A \ox B)\Big) \ox FC \\
          &\cong \int^{X,A,B,C} (\A(X \ox C,-) \times \A(A \ox B,X)) \cdot
                            (F(A \ox B) \ox FC) \\
          &\cong \int^{A,B,C} \Big(\int^X \A(X \ox C,-) \times \A(A \ox
                B,X)\Big) \cdot (F(A \ox B) \ox FC) \\
          &\cong \int^{A,B,C} \A(A \ox B \ox C,-) \cdot F(A \ox B) \ox FC,
                & \text{(Yoneda)}
\end{align*}
and similarly,
\[
    F * F \cong \int^{A,B,C} \A(A \ox B \ox C,-) \cdot FA \ox F(B \ox C).
\]

\begin{proof}[Proof of Theorem~\ref{thm-fun}]
Using the isomorphisms of equation~\eqref{3to1} and Lemma~\ref{lem-2to1} one
of the Frobenius equations may be written as 
\[
\scalebox{0.9}{
    \xygraph{{\dint^{A,B,C} \hspace{-3ex} \A(A \ox B \ox C,-) \cdot
                    F(A \ox B) \ox FC}="1"
        [r(7)] {\dint^{A,B,C} \hspace{-3ex} \A(A \ox B \ox C,-) \cdot
                    FA \ox FB \ox FC}="2"
        "1"[d(2)] {\dint^{A,B,C} \hspace{-3ex} \A(A \ox B \ox C,-) \cdot
                    F(A \ox B \ox C)}="3"
        "2"[d(2)] {\dint^{A,B,C} \hspace{-3ex} \A(A \ox B \ox C,-) \cdot
                    FA \ox F(B \ox C).}="4"
        "1":"2" ^-{\tint 1 \ox i \ox 1}
        "3":"4" ^-{\tint 1 \ox i}
        "1":"3" _-{\tint 1 \ox r}
        "2":"4" ^-{\tint 1 \ox 1 \ox r}}}
\]
This diagram is seen to commute as $F$ is a Frobenius monoidal functor. The
other Frobenius equation follows from a similar diagram.

To prove the second part of the theorem, assume that $\A$ is a closed
monoidal category and that equation~\eqref{eq-asump} holds. The following
calculation verifies the claim.
\begin{align*}
\int^{A,B,C} \A(A  & \ox B \ox C,-) \cdot F(A \ox B \ox C) \\
    & \cong \int^{A,B,C} \A(C,[A \ox B,-]) \cdot F(A \ox B \ox C)
        & \text{($\A$ closed)} \\
    & \cong \int^{A,B} F(A \ox B \ox [A \ox B,-]) & \text{(Yoneda)} \\
    & \cong \int^{X,A,B} \A(X,A \ox B \ox [A \ox B,-]) \cdot FX
        & \text{(Yoneda)} \\
    & \cong \int^{X} \A(X,-) \ox FX & \eqref{eq-asump} \\
    & \cong F & \text{(Yoneda)}
\end{align*}
\end{proof}

%=========================================================================%

%=========================================================================%

\end{document}